\documentclass[11pt,a4paper]{article}
\usepackage[T1]{fontenc}
\usepackage{lmodern}
\usepackage[margin=28mm]{geometry}
\usepackage{amsmath,amssymb,amsthm}
\usepackage{mathrsfs}
\usepackage{microtype}
\usepackage{xurl}
\usepackage[hidelinks]{hyperref}
\usepackage{graphicx}
\numberwithin{equation}{section}

\newtheorem{theorem}{Theorem}
\newtheorem{conjecture}{Conjecture}
\newtheorem{lemma}[theorem]{Lemma}
\newtheorem{proposition}[theorem]{Proposition}
\newtheorem{corollary}[theorem]{Corollary}
\theoremstyle{definition}
\newtheorem{definition}{Definition}
\theoremstyle{remark}

\DeclareMathOperator{\diff}{diff}
\DeclareMathOperator{\AT}{AT}

\title{List-avoiding orientation of graphs and spanning Alon--Tarsi  subdigraphs}
\author{Chenglong Deng$^1$ \and Xuding Zhu$^1$\thanks{ This research is partially supported by Natural Science Foundation of China through Grant:  NSFC 12371359.} }
\date{%
    $^1$School of Mathematical Sciences, Zhejiang Normal University\\[2ex]%
    \today
}

\begin{document}
\maketitle

\begin{abstract}
Given a graph $G$ and a list assignment $F$ 
that assigns to each vertex $v$ a set $F(v)$ of integers, an  orientation $D$ of $G$ is $F$-avoiding if 
 $d_D^+(v) \notin F(v)$ for each vertex $v$.
Akbari, Dalirrooyfard, Ehsani,
Ozeki and Sherkati 
proposed the following conjectures:
(1)  If $G$ has an orientation $D$ such that $|F(v)| \le d_D^+(v)-1$, 
then $G$ has an $F$-avoiding orientation.  (2) 
If $|F(v)| \le \frac{d_G(v)-1}{2}$, then $G$ has an $F$-avoiding orientation. An application of Combinatorial Nullstellensatz relates the list avoiding orientation problem of a graph $G$ to Alon-Tarsi orientation of spanning subgraphs of $G$.
This paper proves that
every finite  loopless digraph $D$ has a   subset $S$ of arcs such that 
the underlying graph $\underline{S}$ of $S$ is a forest, and $D-S$ is an AT-digraph.
Moreover,  every  
orientation $S'$ of $\underline{S}$ extends to an
orientation $D'$ of $\underline{D}$ such that $D$ and $D'$ have the same out-degree sequence,
and $D'-S'$ is an AT-digraph. 
As a consequence, we confirm Conjecture (1), and a slightly weaker version of Conjecture (2): 
 If $|F(v)| \le d_G(v)/2-1$ for each vertex $v$, 
then $G$ has an $F$-avoiding orientation. We also prove the upper bound $o_\ell(G)\le o(G)+2$ conjectured by Abe and Ozeki for the list-orientation number.
\end{abstract}

\section{Introduction and notation}

All graphs and digraphs in this paper are finite and loopless. Parallel edges and arcs are allowed.
Each edge and arc has a distinct label, which is preserved when the arc
is reversed. We write $e=u\to v$ or $e=(u,v)$ to indicate that
the arc $e$ is directed from $u$ to $v$.
The \emph{underlying graph} $\underline D$ of a digraph $D$
is the multigraph obtained by forgetting the directions of the
arcs and retaining one edge with the same label for each arc.
We call $D$ an \emph{orientation} of $\underline D$.
For $v\in V(D)$, let
$A_D^+(v)=\{e:e=v\to u\text{ for some }u\}$,
$d_D^+(v)=|A_D^+(v)|$,
$A_D^-(v)=\{e:e=u\to v\text{ for some }u\}$, and
$d_D^-(v)=|A_D^-(v)|$.

For a digraph $D$, we use $V(D), A(D)$ to denotes its vertex set and arc set, respectively. For simplicity, we may also use $D$ to denote its arc set. For a subset $S\subseteq A(D)$, we also use $S$ to denote the spanning
subdigraph of $D$ with vertex set $V(D)$ and arc set $S$.
We call $S$ a cycle, a path, or a forest in $D$ if
$\underline S$ is a cycle, a path, or a forest, respectively.
Two parallel edges of the underlying graph form a cycle of length two.

\begin{definition}
    A subset $B$ of $A(D)$
is \emph{Eulerian} if $d_B^+(v)=d_B^-(v)$ for every vertex $v$.   
Let $\mathcal E(D)$ be the family of all Eulerian
arc sets of $D$. The {\em difference} of $D$ is defined as 
\[
 \diff(D)=\sum_{B\in\mathcal E(D)}(-1)^{|B|}.
\]
\end{definition}

\begin{definition}\label{def:AT}
A digraph $D$ is an \emph{Alon--Tarsi digraph} (an
\emph{AT-digraph} for short) if $\diff(D)\ne0$. 
For a graph $G$ and a mapping $f: V(G) \to \mathbb N$, we say $G$ is {\em $f$-AT}
if $G$ has an AT-orientation $D$ with $d_D^+(v) < f(v)$ for each vertex $v$. The
\emph{Alon-Tarsi number} of a graph $G$ is
\[
 \AT(G)=\min\{k:\ G\text{ is $k$-AT }\}.
\]
\end{definition}

The concept of Alon-Tarsi number of a graph is motivated by the applications of Combinatorial Nullstellensatz to 
the study of list colouring of graphs. 
It was proved by Alon and Tarsi ~\cite{AT}   that if $G$ is $f$-AT, then $G$ is $f$-choosable. 
This result was extended by Schauz ~\cite{Schauz} who proved that if $G$ is $f$-AT, then $G$ is $f$-paintable.  
It follows that for any graph $G$, 
${\rm ch}(G)\le \chi_P(G) \le \AT(G)$. (The choice number ${\rm ch}(G)$ and the {\em paint number} $\chi_P(G)$ of $G$
are not needed in this paper, and readers are referred to \cite{ZB} for their definitions). 
The Alon-Tarsi number of a graph has become an important graph invariant of independent interests and has been studied extensively in the literature (cf. \cite{ZB}).

If a digraph $D$ is not an AT-digraph, then a natural question is how far is it from being an AT-digraph.
An \emph{AT-obstruction set} of $D$ is a set
$S\subseteq A(D)$ such that $\diff(D-S)\ne0$. Define
\[
 \kappa(D)=\min\{|S|:\ S\subseteq A(D),\ \diff(D-S)\ne0\}.
\]
This minimum exists as $\diff(\emptyset) = 1$. An obstruction set $S$ is \emph{minimum} if  
$|S|=\kappa(D)$, and $S$ is \emph{inclusion-minimal} if no proper subset of $S$
is an obstruction set.   

The following are two main results of this paper.

\begin{theorem}\label{thm:forest}
Every digraph $D$ has a minimum AT-obstruction set that is a forest.
\end{theorem}

\begin{theorem}\label{thm:sameoutdegree}
Let $S$ be an inclusion-minimal AT-obstruction set of $D$.
Every orientation $S'$ of $\underline S$ extends to an orientation
$D'$ of $\underline D$ such that
\[
 d_{D'}^+(v)=d_D^+(v)\qquad(v\in V(D)),
\]
$S'$ is an inclusion-minimal AT-obstruction set of $D'$, and
\[
 |\diff(D'-S')|=|\diff(D-S)|.
\]
\end{theorem}

It was conjectured by Grytczuk \cite{Grytczuk} that any digraph $D$ has an AT-obstruction set $S$ such that   $d_{D-S}^+(v) \ge d_D^+(v)-1$ for each vertex $v$. The combination of Theorems \ref{thm:forest} and \ref{thm:sameoutdegree} yields the following result.

\begin{corollary}
    \label{cor:Grytczuk}
   For every digraph $D$, there is an orientation $D'$ of
$\underline D$ with $d_{D'}^+(v)=d_D^+(v)$ for every
$v\in V(D)$ such that $D'$ has an AT-obstruction set $S$
satisfying $d_{D'-S}^+(v)\ge d_D^+(v)-1$ for every $v\in V(D)$.
\end{corollary}

The spanning AT-subdigraph of a digraph is also related to the list-avoiding orientation of graphs. 
A \emph{forbidden-list assignment} $F$ of a graph $G$ assigns a finite set of
nonnegative integers to every vertex of  $G$. An orientation $D$ of $G$ 
is \emph{$F$-avoiding} if  $d_D^+(v) \notin F(v)$ for every vertex $v$. 
The concept of $F$-avoiding orientation of $G$ was introduced by 
  Akbari, Dalirrooyfard, Ehsani, Ozeki and
Sherkati  in ~\cite{ADEOS}. It was proved in~\cite{ADEOS}    that if   a spanning subgraph $G'$ of $G$ has an AT-orientation $D'$ with $d_{D'}^+(v) \ge |F(v)|$, then $G$ has an $F$-avoiding orientation. The following result is a direct consequence of  Theorems \ref{thm:forest} and \ref{thm:sameoutdegree}:

\begin{corollary}
    \label{cor:Ddegree}
    If $G$ has an orientation $D$ such that 
    $|F(v)|\le d_D^+(v)-1$ for every vertex, then $G$ has an
$F$-avoiding orientation.
\end{corollary}

Corollary \ref{cor:Ddegree} confirms a conjecture proposed in~\cite{ADEOS}, where another conjecture   was also proposed:

\begin{conjecture}
    \label{conj:degree}
    Assume $F$ is a forbidden-list assignment of $G$. If   
    $|F(v)|\le \frac{d_G(v)-1}{2}$ for every vertex, then $G$ has an
$F$-avoiding orientation.
\end{conjecture}

 This conjecture has been studied in a few papers. Akbari et al.~\cite{ADEOS} proved that $G$ has an $F$-avoiding
orientation if $|F(v)|\le d_G(v)/4$ for every vertex $v$.
Bradshaw et al.~\cite{BCMMW} proved that the same conclusion holds
if $|F(v)|\le d_G(v)/3-1$ for every vertex $v$.
Wu, Yang and Zhu~\cite{WYZ} further improved the bound to
$|F(v)|\le d_G(v)/3$ for every vertex $v$.

As a consequence of Corollary \ref{cor:Ddegree},    we prove
a slightly weaker version of Conjecture \ref{conj:degree}. 

\begin{corollary}\label{cor:degree}
If $|F(v)| \le \max\left\{0, \frac{d_G(v)}2 -1\right\}$, then $G$ has an $F$-avoiding orientation. 
\end{corollary}

Following Abe and Ozeki~\cite{AO}, let $o(G)$ be the least
integer $k$ such that $G$ has an orientation of maximum
out-degree at most $k$. The list-orientation number $o_\ell(G)$
is the least positive integer $k$ such that, for every assignment
of lists $L(v)\subseteq\{0,\ldots,d_G(v)\}$ with
$|L(v)|\ge\min\{d_G(v)+1,k\}$ at every vertex, $G$ has an
orientation $D$ satisfying $d_D^+(v)\in L(v)$ at every vertex.

The following corollary confirms a conjecture proposed in 
\cite[Conjecture~1.8]{AO}.

\begin{corollary}\label{cor:list-orientation}
For every finite simple graph $G$, $o_\ell(G)\le o(G)+2$.
\end{corollary}
 
For a path or cycle $C$, choose a traversal. Let $C^+$ be the arcs directed along the traversal
and $C^-$ those directed against it.  Put
\[
 \sigma_C(h)=
 \begin{cases}
  1,&h\in C^+,\\
 -1,&h\in C^-,\\
  0,&h\notin C.
 \end{cases}
\]
A cycle $C$ in $D$ is a {\em directed cycle} if $C^+ = \emptyset$ or $C^- = \emptyset$.
It is well-known that if $B$ is an Eulerian subset of $D$, then $B$ decomposes into the disjoint union of directed cycles.

 For a vertex set $A\subseteq V(D)$,
let $\partial_D^+(A)$ and $\partial_D^-(A)$ be the set of arcs leaving and entering $A$, respectively. 
Define the {\em cut-sign} $\sigma_A(h)$ of $h \in A(D)$ as 
\[
 \sigma_A(h)=
 \begin{cases}
  1,&h\in\partial_D^+(A),\\
 -1,&h\in\partial_D^-(A),\\
  0,&\text{otherwise.}
 \end{cases}
\]

\section{Prescribing directions on an obstruction set}

For an arc $e=u\to v$, let $e^{-1}=v\to u$. If $C$ is a directed cycle
in $D$, write $D^C$ for the orientation obtained by reversing
every arc of $C$. A directed cycle has one incoming and one
outgoing arc at each of its vertices.

\begin{lemma}[Cycle reversal]\label{lem:cycle}
If $C$ is a directed cycle in $D$, then
\[
 \diff(D^C)=(-1)^{|C|}\diff(D).
\]
Moreover, for every $e\in C$,
\begin{equation}\label{eq:cycle-reversal}
 \diff(D^C-e^{-1})
 =(-1)^{|C|}\bigl(\diff(D)-\diff(D-e)\bigr).
\end{equation}
\end{lemma}

\begin{proof}
For $B\in\mathcal E(D)$ define an arc set of $D^C$ by
\[
 \Psi(B)=(B\setminus C)\cup
             \{h^{-1}:h\in C\setminus B\}.
\]
Then  $\Psi $ is a  
bijection from $\mathcal E(D)$ to $\mathcal E(D^C)$.
As $  
 (-1)^{|\Psi(B)|}=(-1)^{|C|}(-1)^{|B|}$, by summing over all $B$ proves the first identity.

For a fixed $e\in C$, the reversed arc $e^{-1}$ is absent from
$\Psi(B)$ exactly when $e\in B$. Therefore the same bijection
restricts to Eulerian arc sets of $D^C-e^{-1}$ and gives
\[
 \begin{aligned}
 \diff(D^C-e^{-1})
 &=(-1)^{|C|}
       \sum_{\substack{B\in\mathcal E(D)\\e\in B}}(-1)^{|B|}\\
 &=(-1)^{|C|}\bigl(\diff(D)-\diff(D-e)\bigr).
 \end{aligned}
\]
\end{proof}

\begin{lemma}\label{lem:one-edge}
Let $S$ be an inclusion-minimal AT-obstruction set of $D$, and
let $e\in S$. There is a directed cycle $C\subseteq D-S+e$
containing $e$ such that, after reversing $C$, the same labelled
edge set $S$ is an inclusion-minimal AT-obstruction set and
\[
 |\diff(D^C-S)|=|\diff(D-S)|.
\]
\end{lemma}

\begin{proof}
Put $H=D-S$. Since $S$ is an obstruction set, $\diff(H)\ne0$.
Since $S\setminus\{e\}$ is a proper subset of $S$, its
inclusion-minimality gives $\diff(H+e)=0$. Partitioning the
Eulerian sets of $H+e$ according to whether they contain $e$,
we obtain
\[
 \sum_{\substack{B\in\mathcal E(H+e)\\e\in B}}(-1)^{|B|}
 =\diff(H+e)-\diff(H)=-\diff(H)\ne0.
\]
Hence $D-S+e$ has a directed cycle $C$
containing $e$. 
Apply~\eqref{eq:cycle-reversal} to $H+e$. The left side is
$\diff(D^C-S)$, because all the other labels of $S$ were already
absent. Consequently,
\[
 \diff(D^C-S)
 =(-1)^{|C|}\bigl(0-\diff(H)\bigr)
 =(-1)^{|C|+1}\diff(H),
\]
which proves the asserted equality and shows that $S$ remains
an obstruction set.

It remains to verify inclusion-minimality. Fix $R\subsetneq S$.
If $e\notin R$, all arcs of $C$ remain in $D-R$, and hence
\[
 \diff(D^C-R)=(-1)^{|C|}\diff(D-R)=0
\]
by Lemma~\ref{lem:cycle}. If $e\in R$, put
$K=D-(R\setminus\{e\})$. The cycle $C$ is contained in $K$.
Both $R\setminus\{e\}$ and $R$ are proper subsets of $S$, so
$\diff(K)=\diff(K-e)=0$. Formula~\eqref{eq:cycle-reversal}
now gives $\diff(D^C-R)=0$. Thus every proper subset of $S$
still fails to be an obstruction set.
\end{proof}

\begin{proof}[Proof of Theorem~\ref{thm:sameoutdegree}]
If $S=\varnothing$, take $D'=D$. Otherwise, whenever an edge of
$S$ has the wrong direction relative to $S'$, apply
Lemma~\ref{lem:one-edge} to that edge in the current orientation.
The directed cycle supplied by the lemma contains no other edge
of $S$. Its reversal therefore corrects the chosen direction and
does not change any direction already prescribed on $S$. The reversal of the directions of all arcs in a directed cycle does not change the out-degree. 
\end{proof}

\section{Cut and cycle identities}

\begin{lemma}[Cut identity]\label{lem:cut}
For every digraph $K$ and every $A\subseteq V(K)$,
\begin{equation}\label{eq:cut}
 \sum_{h\in A(K)}\sigma_A(h)\diff(K-h)
 =\bigl(|\partial_K^+(A)|-|\partial_K^-(A)|\bigr)\diff(K).
\end{equation}
\end{lemma}

\begin{proof}
An Eulerian arc set $B$ enters and leaves $A$ the same number of times. Hence   
\[
 \sum_{h\in B}\sigma_A(h)=0.
\]
An Eulerian set $B$ contributes to $\diff(K-h)$ precisely when
$h\notin B$. Interchanging two finite sums, and then using the
displayed equality, gives
\[
 \begin{aligned}
 \sum_{h\in A(K)}\sigma_A(h)\diff(K-h)
 &=\sum_{B\in\mathcal E(K)}(-1)^{|B|}
                         \sum_{h\in A(K)\setminus B}\sigma_A(h)\\
 &=\sum_{B\in\mathcal E(K)}(-1)^{|B|}
                         \sum_{h\in A(K)}\sigma_A(h)\\
 &=\bigl(|\partial_K^+(A)|-|\partial_K^-(A)|\bigr)
                         \diff(K).
 \end{aligned}
\]
\end{proof}

Let $B$ and $U$ be disjoint arc sets in a fixed
digraph,  and define
\[
 \Phi_B(U)=
 \sum_{\substack{Q\subseteq B\\Q\cup U\text{ is Eulerian}}}
     (-1)^{|Q|}.
\]

\begin{lemma}[Restoration identities]\label{lem:restoration}
For disjoint arc sets $B$ and $U$,
\begin{equation}\label{eq:restoration}
 \begin{aligned}
 \diff(B+U)&=\sum_{L\subseteq U}(-1)^{|L|}\Phi_B(L),\\
 \Phi_B(U)&=\sum_{L\subseteq U}(-1)^{|L|}\diff(B+L).
 \end{aligned}
\end{equation}
\end{lemma}

\begin{proof}
Every Eulerian arc set $E\subseteq B\cup U$ has a unique
decomposition $E=Q\cup L$, where $Q=E\cap B\subseteq B$ and
$L=E\cap U\subseteq U$. The two parts are disjoint, so
$(-1)^{|E|}=(-1)^{|Q|}(-1)^{|L|}$. For a fixed $L$, the sum of
$(-1)^{|Q|}$ over the possible $Q$ is precisely $\Phi_B(L)$.
Summing over all $L\subseteq U$ proves the first identity.

For the second identity, substitute the first into the
right side and collect the coefficient of each $\Phi_B(M)$:
\[
 \begin{aligned}
 \sum_{L\subseteq U}(-1)^{|L|}\diff(B+L)
 &=\sum_{L\subseteq U}(-1)^{|L|}
                   \sum_{M\subseteq L}(-1)^{|M|}\Phi_B(M)\\
 &=\sum_{M\subseteq U}\Phi_B(M)
       (-1)^{|M|}\sum_{M\subseteq L\subseteq U}(-1)^{|L|}.
 \end{aligned}
\]
For a fixed $M$, write $L=M\cup N$ with
$N\subseteq U\setminus M$. Its coefficient becomes
\[
 (-1)^{|M|}(-1)^{|M|}
       \sum_{N\subseteq U\setminus M}(-1)^{|N|}
 =\sum_{N\subseteq U\setminus M}(-1)^{|N|}.
\]
If $M=U$, this sum is $1$, contributed by $N=\varnothing$.
Otherwise, choose an element of $U\setminus M$ and pair each
subset not containing it with the subset obtained by adding it;
the two signs cancel. The coefficient is then $0$. Hence only
the term $M=U$ remains, proving the second identity.
\end{proof}

\begin{lemma}[Cycle identity]\label{lem:values}
Let $k=\kappa(D)$, let $W\subseteq A(D)$ have size $k+1$, and
put $B=D-W$. For every cycle $C\subseteq W$,
\begin{equation}\label{eq:cycle-values}
 \sum_{h\in C}\sigma_C(h)\diff(B+h)=0.
\end{equation}
\end{lemma}

\begin{proof}
Choose the traversal used to define $C^+$ and $C^-$. For every $Q\subseteq B$,
$C^+ \cup Q \in \mathcal{E} (D)$ if and only if $C^- \cup Q \in \mathcal{E}(D)$. Hence 
$\Phi_B(C^+)=\Phi_B(C^-)$. By \eqref{eq:restoration}, 
\[
 \sum_{L\subseteq C^+}(-1)^{|L|}\diff(B+L)
 =\sum_{L\subseteq C^-}(-1)^{|L|}\diff(B+L).
\]
The empty-set term is $\diff(B)$ on both sides and cancels.
For any term with $|L|\ge2$,  
$\diff(D-(W\setminus L))=0$, as the number of deleted arcs is at most $k-1$. 
The only remaining terms are the singleton terms, each with
sign $-1$. Hence,
\[
  \sum_{h\in C^+}\diff(B+h)
 = \sum_{h\in C^-}\diff(B+h).
\]
\end{proof}

\section{Minimum AT-obstruction forests}

\begin{lemma}[Forest exchange]\label{lem:exchange}
Let $k=\kappa(D)$ and let $\varnothing\ne J\subseteq A(D)$.
Let $T\subseteq J$ be a spanning forest of $(V(D),J)$ such that
  $T$ and $J$ have exactly the same connected components.
Suppose an integer $r\ge |T|$ satisfies
\[
 |J|+r-|T|=k.
\]
Define
\[
 \mathscr F=\{F\subseteq A(D): F\text{ is a forest},
                         \ |F|=r,\ T\subseteq F\},
 \qquad H_F=D-(J\cup F).
\]
If $\diff(H_{F_0})\ne0$ for some $F_0\in\mathscr F$, then
there exist $F\in\mathscr F$, $e\in T$, and $f\in A(H_F)$
joining different components of $F$ such that
\begin{equation}\label{eq:exchange}
 \diff(H_F+e-f)\ne0.
\end{equation}
\end{lemma}

\begin{proof}
For any $F\in\mathscr F$, for any $j \in J-T$, $j \notin F$, for otherwise $F$ would contain a cycle. So
 $  F\cap J=T$.
Adding any edge of $J$ to $F$ joins vertices already
connected in $T\subseteq F$, and hence does not merge components.
Hence $J\cup F$ has exactly the same components as $F$. Finally,
\begin{equation}\label{eq:union-size}
 |J\cup F|=|J|+|F|-|J\cap F|=|J|+r-|T|=k.
\end{equation}

Put $a_F=\diff(H_F)$. Suppose, for a contradiction, that the
conclusion fails. Then
\begin{equation}\label{eq:external-zero}
 \diff(H_F+e-f)=0
\end{equation}
for every $F\in\mathscr F$, every $e\in T$, and every
$f\in A(H_F)$ whose endpoints lie in different components of
$F$.

For a fixed $F$, let
\[
 I_F=\{f\in A(H_F):\text{both endpoints of }f
                        \text{ lie in one component of }F\}.
\]
For each $e\in T$, removing $e$ from its tree in $F$ splits
that tree into two smaller trees. Choose either of their vertex
sets and call it $A_{F,e}$. For   $f\in I_F$, let 
\begin{equation}\label{eq:cut-coefficient}
 c_F(e,f)=\sigma_{A_{F,e}}(e)\sigma_{A_{F,e}}(f).
\end{equation}
Note that $c_F(e,f)$ does not depend on the choice of $A_{F,e}$. 
 Indeed, if $f$ lies in the same original $F$-tree as
$e$, exchanging the two pieces negates both cut signs, leaving
their product unchanged. If $f$ lies in a different $F$-tree,
its cut sign is $0$ for either choice.  

By~\eqref{eq:union-size}, the digraph $H_F+e$ is obtained by
deleting exactly $k-1$ arcs from $D$. Thus
$\diff(H_F+e)=0$. Apply Lemma~\ref{lem:cut} to this digraph
and the vertex set $A_{F,e}$. The term from $e$ is
$\sigma_{A_{F,e}}(e)a_F$. Every term from an arc joining
different $F$-components vanishes by~\eqref{eq:external-zero}.
The remaining arcs are exactly those in $I_F$, so
\begin{equation}\label{eq:forest-cut-sum1}
 0=\sigma_{A_{F,e}}(e)a_F
       +\sum_{f\in I_F}\sigma_{A_{F,e}}(f)
                         \diff(H_F+e-f).
\end{equation}
Multiply this equality by $\sigma_{A_{F,e}}(e)$,  and sum over $e\in T$, using $\sigma_{A_{F,e}}(e)^2=1$, gives
\begin{equation}\label{eq:forest-cut-sum}
 0=|T|a_F+\sum_{f\in I_F}\sum_{e\in T}
                       c_F(e,f)\diff(H_F+e-f).
\end{equation}

We multiply~\eqref{eq:forest-cut-sum} by $a_F$ and sum over
$F\in\mathscr F$. We will group the terms indexed by $(F,f)$
according to $W=F+f$. Let
\[
 \mathscr W=\{F\cup\{f\}:F\in\mathscr F,\ f\in I_F\}.
\]
For each such pair, the definition of $I_F$ gives
$f\in A(H_F)=A(D)\setminus(J\cup F)$, and the endpoints of $f$
lie in the same component of $F$. Thus $W=F+f$ has a unique
cycle $C_W$, formed by $f$ and the unique path in $F$ between
its endpoints. Since $F\cap J=T$, we have
\begin{equation}\label{eq:WJ}
 W\cap J=T,\qquad |J\cup W|=|J\cup F|+1=k+1.
\end{equation}

For a fixed $W$, every pair $(F,f)$ yielding $W$ has
$f\in C_W\setminus T$. Conversely, take $h\in C_W\setminus T$
and put $F=W-h$. Then $F$ is a forest with $r$ edges containing
$T$, so $F\in\mathscr F$. Also $h\notin J\cup F$ by
\eqref{eq:WJ}, so $h\in A(H_F)$. The rest of $C_W$ connects the
endpoints of $h$ in $F$, giving $h\in I_F$. Hence the pairs
yielding $W$ are exactly $(W-h,h)$ for $h\in C_W\setminus T$.

For each $W$, fix one traversal of $C_W$ and let
$\sigma_W=\sigma_{C_W}$. Define the signed sum over these pairs by
\[
 b_W=\sum_{h\in C_W\setminus T}\sigma_W(h)a_{W-h}.
\]

Figure~\ref{fig:forest-exchange} illustrates the reindexing
of the pairs $(F,f)$ by $W=F+f$.

\begin{figure}[htbp]
  \centering
  \includegraphics[width=\linewidth]{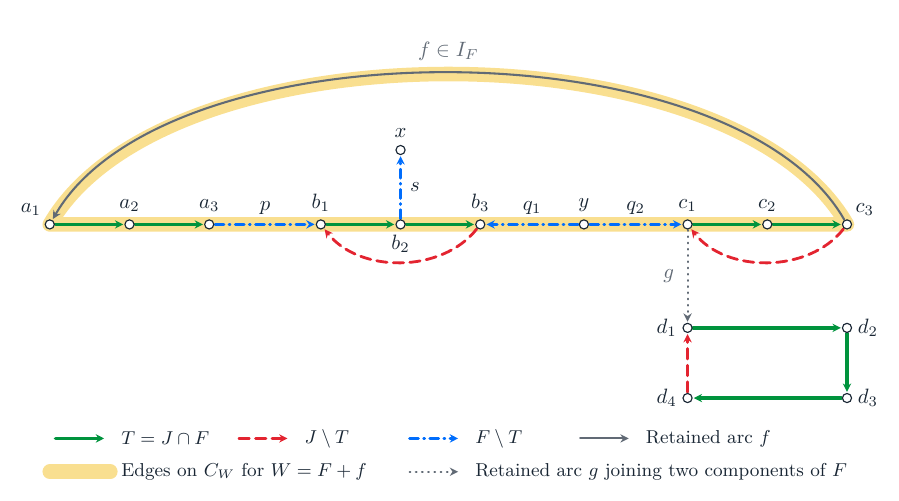}
  \caption{
  The thick solid green arcs form $T=J\cap F$, the dash-dotted
blue arcs form $F\setminus T$, and the dashed red arcs form
$J\setminus T$. Adding the thin solid grey arc $f\in I_F$
to the forest $F$ gives $W$. The highlighted edges form the unique cycle
  $C_W$ of the underlying graph of $W$.
  Here $C_W\setminus T=\{p,q_1,q_2,f\}$, so the four pairs yielding from $W$
  are $(W-p,p)$, $(W-q_1,q_1)$, $(W-q_2,q_2)$, and $(W-f,f)$.
  In the traversal along the $F$-path from $a_1$ to $c_3$ and back
  along $f$, only $q_1$ is traversed against its orientation.
  Thus $b_W=a_{W-p}-a_{W-q_1}+a_{W-q_2}+a_{W-f}$.}
  \label{fig:forest-exchange}
\end{figure}

Put $B_W=D-(J\cup W)$. For $h\in C_W\setminus T$, we have
$h\notin J$, so
\[
 B_W+h=D-\bigl(J\cup(W-h)\bigr)=H_{W-h}.
\]
Apply Lemma~\ref{lem:values} with deleted arc set $J\cup W$ and cycle $C_W\subseteq W$. The sum of the terms indexed by
$C_W\setminus T$ is $b_W$, so
\begin{equation}\label{eq:forest-cycle-sum}
 \sum_{e\in C_W\cap T}\sigma_W(e)\diff(B_W+e)=-b_W.
\end{equation}

Now fix $F\in\mathscr F$ and $f\in I_F$, and put $W=F+f$.
For $e\in T$, the only edge of $F$ crossing the cut of $A_{F,e}$
is $e$. For $h\in C_W$, the product
$\sigma_W(h)\sigma_{A_{F,e}}(h)$ is the cut sign of $h$
when taken in the traversal direction. The traversal enters
and leaves $A_{F,e}$ equally often, so
\[
 0=\sum_{h\in C_W}\sigma_W(h)\sigma_{A_{F,e}}(h)
   =\sigma_W(e)\sigma_{A_{F,e}}(e)
     +\sigma_W(f)\sigma_{A_{F,e}}(f),
\]
where $\sigma_W(e)=0$ when $e\notin C_W$. Since $e$ crosses
the cut and $f$ lies on $C_W$, we have
$\sigma_{A_{F,e}}(e)^2=\sigma_W(f)^2=1$. Multiplying by
$\sigma_{A_{F,e}}(e)\sigma_W(f)$ yields
$c_F(e,f)=-\sigma_W(e)\sigma_W(f)$.
As $H_F-f=B_W$, equation~\eqref{eq:forest-cycle-sum} now gives
\begin{equation}\label{eq:inner-sum}
\begin{aligned}
 \sum_{e\in T}c_F(e,f)\diff(H_F+e-f)
 &=-\sigma_W(f)\sum_{e\in C_W\cap T}
       \sigma_W(e)\diff(B_W+e)\\
 &=\sigma_W(f)b_W.
\end{aligned}
\end{equation}

In the sum of \eqref{eq:forest-cut-sum} over $F\in\mathscr F$,
with weight $a_F$, the contribution indexed by $(F,f)$ is
$\sigma_W(f)a_Fb_W$ by \eqref{eq:inner-sum}.
Grouping these contributions by $W$ and writing $(F,f)=(W-h,h)$ gives
\begin{equation}\label{eq:forest-squares}
\begin{aligned}
 0
 &=|T|\sum_{F\in\mathscr F}a_F^2
   +\sum_{W\in\mathscr W}b_W
       \sum_{h\in C_W\setminus T}\sigma_W(h)a_{W-h}\\
 &=|T|\sum_{F\in\mathscr F}a_F^2
   +\sum_{W\in\mathscr W}b_W^2.
\end{aligned}
\end{equation}
Since $J$ is nonempty and $D$ is loopless, a spanning forest of
$(V(D),J)$ has at least one edge. Thus $|T|>0$. As
$a_{F_0}\ne0$, the first term on the right-hand side is strictly
positive, a contradiction.
\end{proof}

\begin{proof}[Proof of Theorem~\ref{thm:forest}]
For $S\subseteq A(D)$, let $c(S)$ be the number of connected
components of $\underline S$, including isolated vertices.
The rank of $S$ is $r(S)=|V(D)|-c(S)$, which is the number
of edges in any spanning forest of $\underline S$.
Put $k=\kappa(D)$. Among all minimum AT-obstruction sets,
choose $S_0$ with maximum rank. Such a choice exists because
there are finitely many arc sets. If $S_0$ is a forest, there
is nothing to prove; this includes the case $k=0$.

Suppose instead that $S_0$ contains a cycle. Let $J$ be the
set of its non-cut arcs. Choose a spanning tree in each
component of $J$, and let $T$ be the union of their arc sets. Thus $T\subseteq J$ is a
spanning forest and has exactly the same components as $J$.
Put $r=r(S_0)$. The set $J$ is nonempty. Every
edge of $J$ lies on a cycle of $S_0$, and every arc of that
cycle is  a non-cut arc. Hence every arc of $J$ lies
on a cycle contained entirely in $J$.

Set
\[
 F_0=T\cup(S_0\setminus J).
\]
Every edge of $S_0\setminus J$ is a bridge of $\underline{S_0}$.
Thus any cycle of $F_0$ would lie entirely in $T$, which is a forest. Moreover, $F_0$ has the same connected components as $S_0$, because each
edge omitted from $S_0$ belongs to $J\setminus T$ and has its
endpoints connected in $T$. Thus $F_0$ is a spanning forest
of $S_0$ with the additional property that their components
coincide. This property gives $|F_0|=r(S_0)$, and hence
\[
 |F_0|=r=|T|+|S_0\setminus J|,
 \qquad |J|+r-|T|=|S_0|=k.
\]
Moreover, $J\cup F_0=S_0$, so
$\diff(D-(J\cup F_0))=\diff(D-S_0)\ne0$.

Apply Lemma~\ref{lem:exchange} with these choices. It supplies
$F$, $e\in T$, and $f\in A(H_F)$ such that $f$ joins
different components of $F$ and
$\diff(H_F+e-f)\ne0$. Define
\[
 S'=(J\cup F)\setminus\{e\}\cup\{f\}.
\]
Here $e\in J\cup F$ and $f\notin J\cup F$, so $|S'|=k$.
Also $D-S'=H_F+e-f$ is AT. Therefore $S'$ is a minimum
AT-obstruction set.

The edge $e\in T\subseteq J$ lies on a cycle of $J$, so
deleting it from $J\cup F$ does not split any component.
As proved in the exchange lemma, $J\cup F$ and $F$ have the
same components. The added edge $f$ joins two of these
components, thereby reducing the component count by one.
Since $F$ is a forest with $r$ edges, this gives
\[
 r(S')=r(J\cup F)+1=r(F)+1=r+1.
\]
This contradicts the maximal rank of $S_0$. Hence $S_0$ is
a forest.
\end{proof}

\begin{corollary}\label{cor:bounded-loss}
Let $D$ be an orientation of $G$. There exist a minimum forest
AT-obstruction set $S$ of $D$ and an orientation $D^*$ of $G$
such that $d_{D^*}^+(v)=d_D^+(v)$ for every vertex,
$D^*-S$ is AT, and
\[
 d_{D^*[S]}^+(v)\le1\qquad(v\in V(G)).
\]
In particular, the AT-orientation $H=D^*-S$ of the spanning
subgraph $G-S$ satisfies
\begin{equation}\label{eq:bounded-loss}
 d_H^+(v)\ge\max\{0,d_D^+(v)-1\}\qquad(v\in V(G)).
\end{equation}
\end{corollary}

\begin{proof}
Choose $S$ using Theorem~\ref{thm:forest}. In each tree of
$S$, choose a vertex called its root, and direct every edge
towards the root along the unique path to it. Each non-root
vertex then has exactly one outgoing arc of $S$, and the root
has none. This also describes isolated vertices, whose trees
have no edges.

The set $S$ is minimum and hence inclusion-minimal, so
Theorem~\ref{thm:sameoutdegree} realizes these directions in
an orientation $D^*$ of the full graph without changing any
full outdegree. Deleting $S$ leaves an AT-digraph and gives
\[
 d_H^+(v)=d_{D^*}^+(v)-d_{D^*[S]}^+(v)
         \ge d_D^+(v)-1.
\]
Outdegrees are also nonnegative, which gives the maximum
with $0$ in~\eqref{eq:bounded-loss}.
\end{proof}

\section{List-avoiding orientations}

Proposition \ref{thm:criterion} is stated as
Theorem~1.4 in~\cite{BCMMW}. For completeness of this paper, we   include 
a proof of this proposition.

\begin{proposition}[Spanning-subgraph criterion]\label{thm:criterion}
If a spanning subgraph of $G$ has an AT-orientation $H$ with
$d_H^+(v)\ge |F(v)|$ for every vertex, then $G$ has an
$F$-avoiding orientation.
\end{proposition}

\begin{proof}
Put $t_v=d_H^+(v)$ and $m=|A(H)|$, so $\sum_v t_v=m$.
By adding integers to $F(v)$ if needed, we may assume that $|F(v)| = t_v$.  
Extend $H$ arbitrarily to an orientation $D$ of $G$, and let $q_v = d_D^+(v)-d_H^+(v)$.  
Associate to  each   arc $e$ of $H$ a variable $x_e$. For $x: A(H) \to \{0,1\}$,  let $D_x$ be the digraph obtained from 
$D$ by   reversing the direction of those arcs $e$ with $x_e=0$.
Let
\[
  g_v(x) = q_v+
       \sum_{\substack{e\in A_H^+(v) }}x_e
       +\sum_{\substack{e\in A_H^-(v)}}(1-x_e).
\]
Then $d_{D_x}^+(v) = g_v(x)$ for each $x: A(H) \to \{0,1\}$ and each vertex $v$. 
Let
\[
 P(x)=\prod_{v\in V(G)}\prod_{a\in F(v)}(g_v(x)-a).
\]
Then $P(x)\ne0$ if and only if $D_x$ is an $F$-avoiding orientation of $G$. To apply Combinatorial Nullstellensatz, we 
 compute $[\prod_{e\in A(H)}x_e]P(x)$, which is the coefficient  of the highest degree monomial
$\prod_{e\in A(H)}x_e$ in $P(x)$.    
Let $$\tilde{P}(x) = \prod_{v \in V(G)} (\sum_{\substack{e\in A_H^+(v)}}x_e
       -\sum_{\substack{e\in A_H^-(v)  }}x_e)^{t_v}.$$
Then $[\prod_{e\in A(H)}x_e]P(x) = [\prod_{e\in A(H)}x_e]\tilde{P}(x)$.    
To produce the desired monomial in this product, each edge
variable must be chosen at exactly one of its two endpoints,
and exactly $t_v$ edge variables must be chosen at $v$.
Initially assign every edge to its tail in $H$; this assigns
exactly $t_v$ edges to $v$. Let $B$ be the set of edges that
are instead assigned to their heads. The new number assigned
to $v$ is
\[
 t_v-d_B^+(v)+d_B^-(v).
\]
It equals $t_v$ for all $v$ exactly when $B$ is Eulerian.
Each head assignment contributes a negative sign, giving the
sign $(-1)^{|B|}$. For each such endpoint assignment, the
$t_v$ distinct variables assigned to $v$ can be selected from
the $t_v$ identical linear factors in $t_v!$ orders. Hence
\begin{equation}\label{eq:Boolean-coefficient}
[\prod_{e\in A(H)}x_e]P(x)= \left(\prod_{v\in V(G)}t_v!\right)
       \sum_{B\in\mathcal E(H)}(-1)^{|B|}
 =\left(\prod_{v\in V(G)}t_v!\right)\diff(H)\ne0.
\end{equation}
By the Combinatorial Nullstellensatz \cite{AT}, there exists $x: A(H) \to \{0,1\}$ such that $P(x) \ne 0$, and hence $D_x$ is an $F$-avoiding orientation of $G$. 
\end{proof}

\begin{proof}[Proof of Corollary~\ref{cor:Ddegree}]
Apply Corollary~\ref{cor:bounded-loss}
to the given orientation $D$. Its AT-oriented spanning
subgraph $H$ satisfies
\[
 d_H^+(v)\ge d_D^+(v)-1\ge |F(v)|
\]
at every vertex. Proposition~\ref{thm:criterion} gives the conclusion.
\end{proof}
\begin{proof}[Proof of Corollary~\ref{cor:degree}]
 It is well-known and easy to see that every graph $G$ has an orientation $D$ with
\begin{equation}\label{eq:balanced}
 \left\lfloor\frac{d_G(v)}2\right\rfloor
 \le d_D^+(v)\le
 \left\lceil\frac{d_G(v)}2\right\rceil
 \qquad(v\in V(G)).
\end{equation}
Now apply Corollary~\ref{cor:bounded-loss} to this balanced
orientation $D$. The resulting AT-oriented spanning subgraph
satisfies
\[
 d_H^+(v)\ge
 \max\left\{0,\left\lfloor\frac{d_G(v)}2\right\rfloor-1\right\}
 \ge |F(v)|.
\]
Proposition~\ref{thm:criterion}  gives the conclusion.
\end{proof}

\begin{proof}[Proof of Corollary~\ref{cor:list-orientation}]
Put $k=o(G)$. Reverse all arcs of an orientation of $G$
with maximum out-degree at most $k$ to obtain an orientation
$D$ with $d_D^+(v)\ge d_G(v)-k$ for every vertex $v$.

Let $L(v)\subseteq\{0,\ldots,d_G(v)\}$ be lists satisfying
$|L(v)|\ge\min\{d_G(v)+1,k+2\}$ at every vertex, and put
$F(v)=\{0,\ldots,d_G(v)\}\setminus L(v)$. Then
$$
|F(v)|=d_G(v)+1-|L(v)|
\le \max\{0,d_G(v)-k-1\}.
$$
By Corollary~\ref{cor:bounded-loss}, a spanning subgraph of $G$
has an AT-orientation $H$ such that
$$
d_H^+(v)\ge\max\{0,d_D^+(v)-1\}
\ge\max\{0,d_G(v)-k-1\}\ge |F(v)|.
$$
Proposition~\ref{thm:criterion} gives an $F$-avoiding orientation
of $G$, whose out-degree at each vertex $v$ therefore belongs
to $L(v)$.
\end{proof}

\section*{AI Declaration}

The authors used ChatGPT 6 in this research. All AI-assisted suggestions and arguments were
carefully reviewed, revised, and independently verified by the authors. The authors take full
responsibility for the mathematical correctness and content of the manuscript.

\end{document}